\documentclass[12pt]{amsart}

\usepackage{amssymb,verbatim,graphicx}
\usepackage[mathscr]{eucal}
\usepackage{enumerate}
\usepackage{xspace}
\usepackage[thinlines]{easymat}

\newtheorem{theorem}{Theorem}[section]

\newtheorem{proposition}[theorem]{Proposition}

\newtheorem{theoremvoid}{Theorem} 
 
\newtheorem{definition}[theorem]{Definition}

\theoremstyle{definition}
\newtheorem{example}[theorem]{Example}
\newtheorem{remark}[theorem]{Remark}

\numberwithin{equation}{section}

\def\&{\wedge}

\newcommand{\eps}{\varepsilon}

\newcommand{\e}{\mathbf{e}}
\newcommand{\f}{\mathbf{f}}
\newcommand{\bu}{\mathbf{u}}
\newcommand{\bv}{\mathbf{v}}
\newcommand{\bw}{\mathbf{w}}
\newcommand{\bx}{\mathbf{x}}
\newcommand{\R}{\mathbb{R}}

\newcommand{\bb}{\mathbb}

\begin{document}

\title{Totally quasi-umbilic timelike surfaces in $\mathbb{R}^{1,2}$}

\author{Jeanne N. Clelland}
\address{Department of Mathematics, 395 UCB, University of
Colorado,
Boulder, CO 80309-0395}
\email{Jeanne.Clelland@colorado.edu}

\thanks{This research was supported in part by NSF grant DMS-0908456.}

\begin{abstract}
For a regular surface in Euclidean space $\mathbb{R}^3$, umbilic points are precisely the points where the Gauss and mean curvatures $K$ and $H$ satisfy $H^2=K$; moreover, it is well-known that the only totally umbilic surfaces in $\R^3$ are planes and spheres.  But for timelike surfaces in Minkowski space $\mathbb{R}^{1,2}$, it is possible to have $H^2=K$ at a non-umbilic point; we call such points {\em quasi-umbilic}, and we give a complete classification of totally quasi-umbilic timelike surfaces in $\mathbb{R}^{1,2}$.  
\end{abstract}

\maketitle

\section{Introduction}\label{intro-sec}

For a regular surface $\Sigma$ in Euclidean space $\bb{R}^3$, it is well-known (see, e.g., \cite{doCarmo76}) that the shape operator at any point $\bx \in \Sigma$ is a self-adjoint linear operator $S_{\bx}:T_{\bx}\Sigma \to T_{\bx}\Sigma$ and is therefore diagonalizable over $\bb{R}$.  But for a timelike surface $\Sigma$ in 3-dimensional Minkowski space $\R^{1,2}$ (cf. Definition \ref{R12-def}), this is no longer necessarily true; the shape operator is still self-adjoint, but because the metric on $T_{\bx}\Sigma$ is now indefinite, it can have any of three algebraic types: diagonalizable over $\bb{R}$, diagonalizable over $\bb{C}$ but not $\bb{R}$, or non-diagonalizable over $\bb{C}$ with a single null eigenvector (see, e.g., \cite{Magid05}, \cite{Magid09}, \cite{Milnor83}).


An {\em umbilic point} of a regular surface $\Sigma$ in either Euclidean or Minkowski space is a point $\bx \in \Sigma$ for which the second fundamental form $\text{II}_{\bx}$ of $\Sigma$ is a scalar multiple of the first fundamental form $\text{I}_{\bx}$.  In the Euclidean case, the Gauss curvature $K$ and mean curvature $H$ of $\Sigma$ satisfy $H^2 - K \geq 0$, and the umbilic points are precisely those points where $H^2 - K = 0$.  But for a timelike surface in Minkowski space, the quantity $H^2 - K$ can take on any real value.  Specifically:
\begin{enumerate}
\item If the shape operator is diagonalizable over $\bb{R}$, then $H^2 - K \geq 0$, with $H^2-K$=0 precisely at umbilic points.
\item If the shape operator is diagonalizable over $\bb{C}$ but not $\bb{R}$, then $H^2 - K < 0$.
\item If the shape operator is non-diagonalizable over $\bb{C}$, then $H^2-K=0$.
\end{enumerate}
Because of this relationship between $H$ and $K$ in Case (3), we make the following definition:

\begin{definition}\label{quasi-umbilic-def}
Let $\Sigma$ be a regular timelike surface in $\bb{R}^{1,2}$.  A point $\bx \in \Sigma$ will be called {\em quasi-umbilic} if the shape operator $S_{\bx}:T_{\bx}\Sigma \to T_{\bx}\Sigma$ is non-diagonalizable over $\bb{C}$.
\end{definition}

\begin{definition}\label{tqu-def}
A regular timelike surface $\Sigma$ in $\bb{R}^{1,2}$ will be called {\em totally quasi-umbilic} if every point $\bx \in \Sigma$ is quasi-umbilic.
\end{definition}

The goal of this paper is to give a classification of totally quasi-umbilic timelike surfaces in $\R^{1,2}$.  Our main result is:

\begin{theoremvoid}[cf. Theorem \ref{main-thm}]
Let $\Sigma$ be a totally quasi-umbilic, regular timelike surface in $\R^{1,2}$. Then $\Sigma$ is a ruled surface whose rulings are all null lines, with the additional property that any null curve $\alpha$ in $\Sigma$ which is transverse to the rulings is nondegenerate (i.e., the vectors $\alpha', \alpha''$ are linearly independent at each point).  Conversely, given any nondegenerate null curve $\alpha(u)$ in $\R^{1,2}$ and any null vector field $\bar{\f}_2(u)$ along $\alpha$ which is linearly independent from $\alpha'(u)$ for all $u$, the immersed ruled surface
\[ \bx(u,v) = \alpha(u) + v \bar{\f}_2(u) \]
is totally quasi-umbilic, possibly containing a curve of umbilic points.  
\end{theoremvoid}

The paper is organized as follows: in \S \ref{Minkowski-review-sec} we briefly review basic notions from the geometry of timelike surfaces in $\R^{1,2}$. In \S \ref{fundamental-forms-sec}, we use Cartan's method of moving frames to derive and solve a system of PDEs whose solutions provide a complete description of the Maurer-Cartan forms for totally quasi-umbilic surfaces in local null coordinates.  In \S \ref{main-result-sec}, we solve the resulting Maurer-Cartan equations to find explicit local parametrizations for all such surfaces and thereby prove Theorem \ref{main-thm}.  Finally, in \S \ref{examples-sec}, we give some examples of totally quasi-umbilic surfaces.

\section{Brief review of geometry of timelike surfaces in $\R^{1,2}$}\label{Minkowski-review-sec}

\begin{definition}\label{R12-def}
Three-dimensional {\em Minkowski space} is the manifold $\R^{1,2}$ defined by
\[ \R^{1,2} = \{\bx = (x^0, x^1, x^2): x^0, x^1, x^2 \in \bb{R} \}, \]
with the indefinite inner product $\langle \cdot, \cdot \rangle$ defined on each tangent space $T_{\bx}\R^{1,2}$ by
\[ \langle \bv, \bw \rangle = v^0 w^0 - v^1 w^1 - v^2 w^2. \]
\end{definition}

\begin{definition}
A nonzero vector $\bv \in T_{\bx}\bb{R}^{1,2}$ is called: 
\begin{itemize}
\item {\em spacelike} if $\langle \bv, \bv \rangle < 0$;
\item {\em timelike} if $\langle \bv, \bv \rangle > 0$;
\item {\em lightlike} or {\em null} if $\langle \bv, \bv \rangle = 0$.
\end{itemize}
\end{definition}

\begin{figure}[h]
\includegraphics[width=2in]{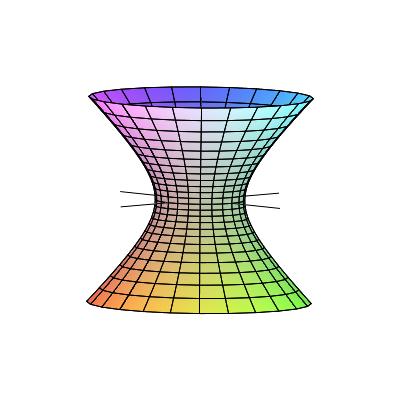}
\includegraphics[width=2in]{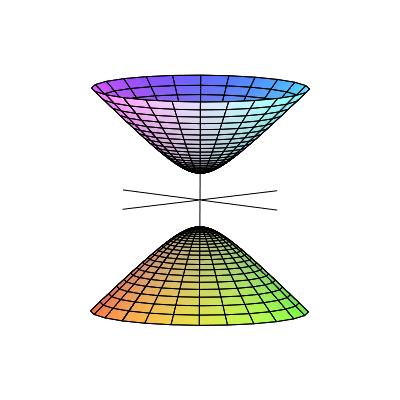}
\includegraphics[width=2in]{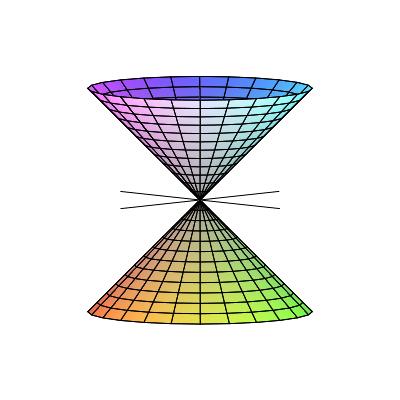}
\caption{``Spheres" in $\R^{1,2}$ with $c < 0, \ c > 0, \ c = 0$}
\label{Minkowski-spheres}
\end{figure}

For a given real number $c$, the ``sphere" consisting of all vectors $\bv$ with $\langle \bv, \bv \rangle = c$ is:
\begin{itemize}
\item a hyperboloid of one sheet if $c<0$;
\item a hyperboloid of two sheets if $c>0$;
\item a cone, called the {\em light cone} or {\em null cone} if $c=0$. \end{itemize}
(See Figure \ref{Minkowski-spheres}; note that the $x^0$-axis is drawn as the vertical axis.)

\begin{definition}
A regular surface $\Sigma \subset \bb{R}^{1,2}$ is called:
\begin{itemize}
\item {\em spacelike} if the restriction of $\langle \cdot, \cdot \rangle$ to each tangent plane $T_{\bx}\Sigma$ is negative definite;
\item {\em timelike} if the restriction of $\langle \cdot, \cdot \rangle$ to each tangent plane $T_{\bx}\Sigma$ is indefinite;
\item {\em lightlike} if the restriction of $\langle \cdot, \cdot \rangle$ to each tangent plane $T_{\bx}\Sigma$ is degenerate.
\end{itemize}
\end{definition}

The Minkowski isometry group consists of all transformations $\varphi: \R^{1,2} \to \R^{1,2}$ of the form
\[ \varphi(\bx) = A\bx + \mathbf{b}, \]
where $A \in O(1,2)$ and $\mathbf{b} \in \R^{1,2}$.  Cartan's method of moving frames can be used to compute local invariants for timelike surfaces under the action of this isometry group.  Given a regular timelike surface $\Sigma \subset \R^{1,2}$, an {\em adapted orthonormal frame} at a point $\bx \in \Sigma$ consists of three mutually orthogonal vectors $(\e_0, \e_1, \e_2)$ in $T_{\bx}\R^{1,2}$ such that:
\begin{itemize}
\item $\langle \e_0, \e_0 \rangle = 1, \ \ \ \langle \e_1, \e_1 \rangle = \langle \e_2, \e_2 \rangle = -1$;
\item $\e_0, \e_1$ span the tangent plane $T_{\bx}\Sigma$.
\end{itemize}
An {\em adapted orthonormal frame field} along $\Sigma$ consists of smooth vector fields $(\e_0, \e_1, \e_2)$ along $\Sigma$ such that for each $\bx \in \Sigma$, the vectors $(\e_0(\bx), \e_1(\bx), \e_2(\bx))$ form an adapted orthonormal frame at $\bx$.
If the surface $\Sigma$ is given an orientation, then the unit normal vector field $\e_2$ along $\Sigma$ is completely determined. In this case, it is convenient to consider only {\em oriented} orthonormal frame fields, defined by the condition that $\e_0 \times \e_1 = -\e_2$, where the Minkowski cross product is given by
\begin{equation}
 \bv \times \bw = (v^1 w^2 - v^2 w^1, \, v^0 w^2 - v^2 w^0, \, v^1 w^0 - v^0 w^1 ) \label{Minkowski-cross-product}
\end{equation}
for $\bv = (v^0, v^1, v^2), \ \bw = (w^0, w^1, w^2)$.  

\begin{remark}
These sign conventions are chosen in order to satisfy the orientation condition
\[ \langle \bu, \bv \times \bw \rangle = \det \begin{bmatrix} \bu & \bv & \bw  \end{bmatrix} \]
for any three vectors $\bu, \bv, \bw$, as well as the condition that the standard basis vectors
\[ \underline{\e}_0 = (1,0,0), \qquad \underline{\e}_1 = (0,1,0), \qquad \underline{\e}_2 = (0,0,1) \]
form an oriented orthonormal frame.
\end{remark}

Given an oriented adapted orthonormal frame field $(\e_0, \e_1, \e_2)$ along $\Sigma$, any other oriented adapted orthonormal frame field $(\tilde{\e}_0, \tilde{\e}_1, \tilde{\e}_2)$ along $\Sigma$ is given by
\begin{equation*}
 \begin{bmatrix} \tilde{\e}_0 & \tilde{\e}_1 & \tilde{\e}_2 
\end{bmatrix} = \begin{bmatrix} \e_0 & \e_1 & \e_2 \end{bmatrix} 
\begin{bmatrix} \cosh \theta & \sinh \theta & 0 \\[0.05in] \sinh \theta & 
\cosh \theta & 0 \\[0.05in] 0 & 0 & 1  \end{bmatrix} 
\end{equation*}
for some real-valued function $\theta$ on $\Sigma$.

As in the Euclidean case, the normal vector field $\e_2$ of an oriented adapted orthonormal frame field along $\Sigma$ defines the {\em Gauss map} from $\Sigma$ to the Minkowski sphere $\langle \bv, \bv \rangle = 1$, and at any point $\bx \in \Sigma$, the differential $d\e_2$ of the Gauss map is a self-adjoint linear operator
\[ d\e_2: T_{\bx}\Sigma \to T_{\bx}\Sigma. \]
The {\em shape operator} of $\Sigma$ at any point $\bx \in \Sigma$ is defined to be 
\[ S_{\bx} = -d\e_2: T_{\bx}\Sigma \to T_{\bx}\Sigma. \]
The first and second fundamental forms $\text{I}_{\bx}, \text{II}_{\bx}:T_{\bx}\Sigma \to \R$ are defined by
\begin{align*}
\text{I}_{\bx}(\bv) & = \langle \bv, \bv \rangle, \\
\text{II}_{\bx}(\bv) & = \langle S_{\bx}(\bv), \bv \rangle,
\end{align*}
and the Gauss curvature $K$ and mean curvature $H$ at $\bx$ are defined to be
\[ K(\bx) = \det \left( S_{\bx} \right), \qquad  H(\bx) =  \tfrac{1}{2} \text{tr} \left( S_{\bx} \right). \]
A point $\bx \in \Sigma$ is called {\em umbilic} if the quadratic form $\text{II}_{\bx}$ is a scalar multiple of $\text{I}_{\bx}$. It is straightforward to show that $\bx$ is umbilic if and only if $S_{\bx}$ is a scalar multiple of the identity; moreover, if $\bx$ is an umbilic point, then $H(\bx)^2 - K(\bx) = 0$.

The shape operator $S_{\bx}$ is a self-adjoint linear operator; i.e., for any vectors $\bv, \bw \in T_{\bx}\Sigma$, we have
\[ \langle S_{\bx}(\bv), \bw \rangle = \langle \bv,  S_{\bx}(\bw) \rangle. \]
This is equivalent to the condition that the matrix representation of $S_{\bx}$ with respect to any orthonormal basis $(\e_0, \e_1)$ for $T_{\bx}\Sigma$ is {\em skew-symmetric}.  So unlike in the Euclidean case, where $S_{\bx}$ is always diagonalizable over $\R$, any of three possibilities may occur: $S_{\bx}$ may be diagonalizable over $\bb{R}$, diagonalizable over $\bb{C}$ but not $\bb{R}$, or non-diagonalizable over $\bb{C}$ with a repeated real eigenvalue and a 1-dimensional eigenspace.

If $S_{\bx}$ is not diagonalizable over $\bb{C}$, then there exists a basis for $T_{\bx}\Sigma$ with respect to which the matrix for $S_{\bx}$ has the form
\[ \begin{bmatrix} \lambda & 1 \\[0.05in] 0 & \lambda \end{bmatrix} \]
for some $\lambda \in \R$.  Thus we have $K(\bx) = \lambda^2, \ H(\bx) = \lambda$, and so $H(\bx)^2 - K(\bx) = 0$.  Because of this similarity with umbilic points, we call such points $\bx$ {\em quasi-umbilic} (cf. Definition \ref{quasi-umbilic-def}).  Moreover, the unique eigenvector for $S_{\bx}$ is necessarily a null vector; to see this, note that with respect to any orthonormal basis $(\e_0, \e_1)$ for $T_{\bx}\Sigma$, the matrix for $S_{\bx}$ is skew-symmetric and so has the form
\[ \begin{bmatrix} a & b \\[0.05in] -b & c \end{bmatrix}. \]
The condition that $S_{\bx}$ has a repeated eigenvalue is
\[ (a-c)^2 = 4 b^2, \]
and $S_{\bx}$ has a 1-dimensional eigenspace (and hence $\bx$ is quasi-umbilic) if and only if $a \neq c$. In this case the repeated eigenvalue is $\lambda = \tfrac{1}{2}(a+c)$, with eigenvector $\e_0 - \e_1$, which is clearly null.

It is well-known that totally umbilic surfaces in Euclidean space $\R^3$ are necessarily contained in either planes or spheres.  The analogous result for timelike surfaces in $\R^{1,2}$ is:

\begin{theorem}
Suppose that a regular timelike surface $\Sigma \subset \R^{1,2}$ is totally umbilic. Then $\Sigma$ is contained in either a plane or a hyperboloid of one sheet (i.e., a Minkowski ``sphere").
\end{theorem}

For the remainder of this paper, we will consider the classification problem for totally quasi-umbilic timelike surfaces (cf. Definition \ref{tqu-def}).

\section{Null frames and Maurer-Cartan forms}\label{fundamental-forms-sec}

We will approach the classification problem for totally quasi-umbilic timelike surfaces via the method of moving frames.  It turns out to be more convenient to use null frames rather than orthonormal frames, and so we make the following definition:

\begin{definition}
An {\em adapted null frame field} along a timelike surface $\Sigma \subset \R^{1,2}$ consists of linearly independent, smooth vector fields $(\f_1, \f_2, \f_3)$ along $\Sigma$ such that for each $\bx \in \Sigma$,
\begin{itemize}
\item $\f_1(\bx), \f_2(\bx)$ are null vectors which span the tangent space $T_{\bx}\Sigma$ at each point $\bx \in \Sigma$;
\item $\langle \f_1(\bx), \f_2(\bx) \rangle = 1$;
\item $\f_3(\bx)$ is orthogonal to $T_{\bx}\Sigma$ at each point $\bx \in \Sigma$, with $\langle \f_3(\bx), \f_3(\bx) \rangle = -1$.
\end{itemize}
\end{definition}
For instance, if $(\e_0, \e_1, \e_2)$ is an adapted orthonormal frame field along $\Sigma$, then the vector fields
\[
\f_1  = \tfrac{1}{\sqrt{2}}(\e_0 + \e_1), \qquad
\f_2  = \tfrac{1}{\sqrt{2}}(\e_0 - \e_1), \qquad
\f_3  = \e_2 
\]
form an adapted null frame field along $\Sigma$. 

Let $\eta^i, \eta^i_j, \ 1 \leq i,j \leq 3$ be the Maurer-Cartan forms on $\Sigma$ associated to an adapted null frame field along $\Sigma$.  These forms are defined by the equations (note: all indices range from 1 to 3, and we use the Einstein summation convention)
\begin{align}
d\bx & = \f_i \eta^i, \label{MC-forms-def}   \\ 
d\f_i & =  \f_j \eta^j_i, \notag
\end{align}
and they satisfy the Maurer-Cartan structure equations
\begin{align}
d\eta^i & = - \eta^i_j \& \eta^j,  \label{structure-eqns} \\
d\eta^i_j & = - \eta^i_k \& \eta^k_j. \notag
\end{align}
(See \cite{IL03} for a discussion of Maurer-Cartan forms and their structure equations.)
Differentiating the inner product relations
\begin{alignat*}{3}
\langle \f_1, \f_1 \rangle & = 0, & \qquad \langle \f_2, \f_2 \rangle & = 0, & \qquad \langle \f_3, \f_3 \rangle & = -1,   \\[0.1in]
\langle \f_1, \f_2 \rangle & = 1, & \qquad \langle \f_1, \f_3 \rangle & = 0, & \qquad \langle \f_2, \f_3 \rangle & = 0 
\end{alignat*}
yields the following relations among the Maurer-Cartan forms:
\begin{gather}
\eta^1_2 = \eta^2_1 = \eta^3_3 = 0, \label{MC-forms-relations}  \\[0.1in]
\eta^2_2 = -\eta^1_1, \qquad \eta^1_3 = \eta^3_2, \qquad \eta^2_3 = \eta^3_1. \notag
\end{gather}

From the equation
\[ d\bx = \f_1 \eta^1 + \f_2 \eta^2 + \f_3 \eta^3 \]
and the fact that $d\bx$ takes values in $T_{\bx}\Sigma$, it follows that $\eta^3 = 0$.  Moreover, $\eta^1$ and $\eta^2$ are linearly independent 1-forms which form a basis for the cotangent space $T^*_{\bx}\Sigma$ at each point $\bx \in \Sigma$.  

Differentiating the equation $\eta^3 = 0$ yields
\[ 0 = d\eta^3 = -\eta^3_1 \& \eta^1 - \eta^3_2 \& \eta^2. \]
By Cartan's Lemma (see \cite{IL03}), there exist functions $k_{ij} = k_{ji}$ on $\Sigma$  such that
\begin{equation}
 \begin{bmatrix} \eta^3_1 \\[0.1in] \eta^3_2 \end{bmatrix} = -\begin{bmatrix} 
k_{11} & k_{12} \\[0.1in] k_{12} & k_{22} \end{bmatrix} \begin{bmatrix} 
\eta^1 \\[0.1in] \eta^2 \end{bmatrix}. \label{kij-def}
\end{equation}
(The minus sign is included for convenience in what follows.)
Compare the matrix in \eqref{kij-def} to that of the shape operator at any point $\bx \in \Sigma$; using \eqref{MC-forms-def}, \eqref{MC-forms-relations}, and \eqref{kij-def}, we have:
\begin{align*}
 S_{\bx} & = -d\f_3 \\
 & = -\f_1 \eta^1_3 - \f_2 \eta^2_3 \\
 & = -\f_1 \eta^3_2 - \f_2 \eta^3_1 \\
 & = \f_1 (k_{12} \eta^1 + k_{22} \eta^2) + \f_2 (k_{11} \eta^1 + k_{12} \eta^2) \\
 & = (k_{12} \f_1 + k_{11} \f_2) \eta^1 + (k_{22} \f_1 + k_{12} \f_2) \eta^2.
\end{align*}
This means that
\[ S_{\bx}(\f_1) = k_{12} \f_1 + k_{11} \f_2, \qquad  S_{\bx}(\f_2) = k_{22} \f_1 + k_{12} \f_2, \]
and so the matrix of $S_{\bx}$ with respect to the basis $(\f_1, \f_2)$ for $T_{\bx}\Sigma$ is
\[ \begin{bmatrix} k_{12} & k_{22} \\[0.1in] k_{11} & k_{12} \end{bmatrix}. \]
Thus we have
\begin{gather*}
 H = k_{12}, \qquad K = k_{12}^2 - k_{11} k_{12}, \\
 H^2 - K = k_{11}k_{22}.
\end{gather*}
It follows that $\Sigma$ is totally quasi-umbilic if and only if $k_{11} k_{22} \equiv 0$ on $\Sigma$ and $k_{11}, k_{22}$ are not both zero at any point $\bx \in \Sigma$.  Without loss of generality, we may assume that $k_{22}=0, \, k_{11} \neq 0$.

Now we will examine how the functions $k_{ij}$ transform if we make a change of adapted null frame field.  So suppose that $(\tilde{\f}_1, \tilde{\f}_2, \tilde{\f}_3)$ is any other adapted null frame field along $\Sigma$.  Then we must have
\begin{equation}
 \begin{bmatrix} \tilde{\f}_1 & \tilde{\f}_2 & \tilde{\f}_3
\end{bmatrix} = \begin{bmatrix} \f_1 & \f_2 & \f_3 \end{bmatrix} 
\begin{bmatrix} \eps_1 e^{\theta} & 0 & 0 \\[0.1in] 0 & 
\eps_1 e^{-\theta} & 0 \\[0.1in] 0 & 0 & \eps_2  \end{bmatrix} \label{frame-change-eqn}
\end{equation}
for some function $\theta$ on $\Sigma$, where $\eps_1, \eps_2 = \pm 1$.  (Note that we cannot exchange $\f_1$ and $\f_2$ due to our assumption that $k_{11} \neq 0, \ k_{22} = 0$.)  Under such a transformation, the Maurer-Cartan forms $\tilde{\eta}^i, \tilde{\eta}^i_j$ associated to the frame field $(\tilde{\f}_1, \tilde{\f}_2, \tilde{\f}_3)$ satisfy the conditions:
\begin{gather*}
\begin{bmatrix} \tilde{\eta}^1 \\[0.1in] \tilde{\eta}^2 \end{bmatrix} = 
\begin{bmatrix} \eps_1 e^{-\theta} & 0 \\[0.1in] 0 & \eps_1 e^{\theta} \end{bmatrix}
\begin{bmatrix} \eta^1 \\[0.1in]  \eta^2 \end{bmatrix} = \begin{bmatrix} \eps_1 e^{-\theta} \eta^1 \\[0.1in]  \eps_1 e^{\theta} \eta^2 \end{bmatrix}, \\[0.2in]
\begin{bmatrix} \tilde{\eta}^1_3 \\[0.1in] \tilde{\eta}^2_3 \end{bmatrix} = 
\eps_2 \begin{bmatrix} \eps_1 e^{-\theta} & 0 \\[0.1in] 0 & \eps_1 e^{\theta} \end{bmatrix}
\begin{bmatrix} \eta^1_3 \\[0.1in]  \eta^2_3 \end{bmatrix} = \begin{bmatrix} \eps_1 \eps_2 e^{-\theta} \eta^1_3 \\[0.1in]  \eps_1 \eps_2 e^{\theta} \eta^2_3 \end{bmatrix} .
\end{gather*}
Therefore,
\[ \begin{bmatrix} \tilde{\eta}^3_1 \\[0.1in]  \tilde{\eta}^3_2 \end{bmatrix} = \begin{bmatrix} \tilde{\eta}^2_3 \\[0.1in]  \tilde{\eta}^1_3 \end{bmatrix} = \begin{bmatrix} \eps_1 \eps_2 e^{\theta} \eta^2_3 \\[0.1in]  \eps_1 \eps_2 e^{-\theta} \eta^1_3  \end{bmatrix} = \begin{bmatrix} \eps_1 \eps_2 e^{\theta} \eta^3_1 \\[0.1in]  \eps_1 \eps_2 e^{-\theta} \eta^3_2  \end{bmatrix}, \]
and the functions $\tilde{k}_{ij}$ associated to the frame field $(\tilde{\f}_1, \tilde{\f}_2, \tilde{\f}_3)$ are given by
\begin{equation*}
\begin{bmatrix} \tilde{k}_{11} & \tilde{k}_{12} \\[0.1in] \tilde{k}_{12} & 
\tilde{k}_{22} \end{bmatrix} = \eps_2 
\begin{bmatrix} \eps_1 e^{\theta} & 0 \\[0.1in] 0 & \eps_1 e^{-\theta} \end{bmatrix} \begin{bmatrix} k_{11} & k_{12} \\[0.1in] k_{12} & k_{22} \end{bmatrix} \begin{bmatrix} \eps_1 e^{\theta} & 0 \\[0.1in] 0 & \eps_1 e^{-\theta} \end{bmatrix} =
\eps_2 \begin{bmatrix} e^{2\theta} k_{11} & k_{12} \\[0.1in] k_{12} & e^{-2\theta} k_{22} \end{bmatrix} . 
\end{equation*}
Since $k_{11} \neq 0$, it follows that there exists an adapted null frame field $(\f_1, \f_2, \f_3)$ along $\Sigma$ with $k_{11}\equiv 1$, and this frame field is unique up to the discrete transformation
$(\tilde{\f}_1, \tilde{\f}_2, \tilde{\f}_3) = (\eps_1 \f_1, \eps_1 \f_2, \f_3)$.  
For such an adapted null frame field, the first and second fundamental forms of $\Sigma$ are:
\begin{gather}
\text{I} = \langle d\bx, d\bx \rangle =  2 \eta^1 \eta^2, \label{fundamental-forms} \\
\text{II} = - \langle d\f_3, d\bx \rangle =  -(\eta^3_1 \eta^1 + \eta^3_2 \eta^2) =
(\eta^1)^2 + 2H \eta^1 \eta^2.  \notag
\end{gather}
(Recall that $H = k_{12}$.)

Now let $u,v$ be local null coordinates on $\Sigma$; note that these are well-defined up to transformations of the form
\[ u \to \tilde{u}(u), \qquad v \to \tilde{v}(v), \]
with $\tilde{u}'(u), \tilde{v}'(v) \neq 0$.  Then there exist functions $f, g$ on $\Sigma$ such that
\begin{equation}
 \eta^1 = e^f\, du, \qquad \eta^2 = e^g\, dv, \label{MC-forms-1}
\end{equation}
and therefore 
\begin{gather}
\text{I} = 2 e^{f+g}\, du\, dv, \label{fundamental-forms-coords} \\
\text{II} = e^{2f}\, du^2 + 2He^{f+g}\, du\, dv. \notag
\end{gather}
From \eqref{kij-def}, we also have
\begin{equation} 
\eta^3_1 = \eta^2_3 = -(e^f\, du + H e^g\, dv), \qquad
\eta^3_2 = \eta^1_3 = -H e^f\, du. \label{MC-forms-2}
\end{equation}
The structure equations \eqref{structure-eqns} for $d\eta^1$ and $d\eta^2$ can be used to show that
\begin{equation}
 \eta^1_1 = -\eta^2_2 = g_u\, du - f_v\, dv, \label{MC-forms-3}
\end{equation}
and according to \eqref{MC-forms-relations}, the remaining Maurer-Cartan forms are zero.

Finally, the remaining structure equations in \eqref{structure-eqns}---which may be interpreted as the Gauss and Codazzi equations for timelike surfaces in $\R^{1,2}$---imply that the functions $f, g, H$ satisfy the following PDE system:
\begin{gather}
H_v = 0, \label{PDE1} \\
H_u = 2f_v e^{f-g}, \label{PDE2} \\
(f+g)_{uv} = H^2 e^{f+g}. \label{PDE3}
\end{gather}
Conversely, any solution to this system gives rise (at least locally) to the Maurer-Cartan forms \eqref{MC-forms-1}, \eqref{MC-forms-2}, \eqref{MC-forms-3} of a totally quasi-umbilic timelike surface.

Equation \eqref{PDE1} implies that $H$ is a function of $u$ alone.  In order to solve the remaining equations, we will divide into cases based on whether $H(u)$ is zero or nonzero. 

\subsection{Case 1: $H(u) \equiv 0$}

Then equations \eqref{PDE2} and \eqref{PDE3} reduce to:
\begin{gather*}
f_v = 0 ,  \\
(f+g)_{uv} = 0. 
\end{gather*}
Thus $f = f(u)$ is a function of $u$ alone, and 
\[ f + g = \phi(u) + \psi(v) \]
for some functions $\phi(u), \psi(v)$.  

The first and second fundamental forms of $\Sigma$ now take the form
\begin{gather*}
\text{I} = 2 e^{\phi(u)} e^{\psi(v)}\, du\, dv, \\
\text{II} = e^{2 f(u)}\, du^2.
\end{gather*}
By making a change of coordinates $(u,v) \to (\tilde{u}, \tilde{v})$ with
\[ d\tilde{u} = e^{\phi(u)}\, du, \qquad d\tilde{v} = e^{\psi(v)}\, dv \]
we can assume that $\phi(u) = \psi(v) = 0$, and so
\begin{gather}
\text{I} = 2\, du\, dv, \label{case1-fundamental-forms} \\
\text{II} = e^{2 f(u)}\, du^2. \notag
\end{gather}
From \eqref{fundamental-forms}, it follows that
\begin{equation}
 \eta^1 = e^{f(u)}\, du, \qquad \eta^2 = e^{-f(u)}\, dv, \label{case1-MC-forms-1}
\end{equation}
and equations \eqref{MC-forms-2}, \eqref{MC-forms-3} become:
\begin{gather} 
\eta^3_1 = \eta^2_3 = -e^f\, du , \qquad
\eta^3_2 = \eta^1_3 = 0, \label{case1-MC-forms-2} \\
 \eta^1_1 = -\eta^2_2 = -f'(u)\, du . \label{case1-MC-forms-3}
\end{gather}
Thus solutions of the PDE system \eqref{PDE1}-\eqref{PDE3} with $H(u) \equiv 0$ depend on one arbitrary function of one variable---i.e., the function $f(u)$.

\subsection{Case 2: $H(u) \neq 0$}

Let $h(u) = \ln |H(u)|$.  Then we can write equation \eqref{PDE3} as
\[ (f+g)_{uv} = e^{f + g + 2 h(u)}, \]
or equivalently,
\[ \left( f+g + 2 h(u) \right)_{uv} = e^{f + g + 2 h(u)}. \]
Thus the function
\[ z(u,v) = f(u,v) + g(u,v) + 2h(u) \]
is a solution to Liouville's equation
\begin{equation}
  z_{uv} = e^z. \label{Liouville's-eqn}
\end{equation}
The general solution to \eqref{Liouville's-eqn} is
\begin{equation} 
z(u,v) = \ln \left( \frac{2 \phi'(u) \psi'(v)}{(\phi(u) + \psi(v))^2} \right), \label{Liouville-gen-soln}
\end{equation}
where $\phi(u), \psi(v)$ are arbitrary functions with $\phi'(u), \psi'(v) \neq 0$; thus we have
\begin{equation}
f(u,v) + g(u,v) + 2h(u) = \ln \left( \frac{2 \phi'(u) \psi'(v)}{(\phi(u) + \psi(v))^2} \right) \label{apply-Liouville-gen-soln}
\end{equation}
for some functions $\phi(u), \psi(v)$ with $\phi'(u), \psi'(v) \neq 0$.

\begin{remark}
The expression \eqref{Liouville-gen-soln} for the general solution to \eqref{Liouville's-eqn} can be derived from the B\"acklund transformation
\begin{equation}\label{Liouvillewave}
\begin{aligned}
z_u - w_u &=  2e^{\frac{(z+w)}{2}}, \\
z_v + w_v &= e^{\frac{(z-w)}{2}},
\end{aligned}
\end{equation}
which relates solutions $z(u,v)$ of \eqref{Liouville's-eqn} to solutions $w(u,v)$ of the wave equation $w_{uv}=0$ (see, e.g., \cite{CI09}).  Substituting the general solution
\[ w(u,v) = \rho(u) + \sigma(v) \]
of the wave equation into \eqref{Liouvillewave} and solving the resulting overdetermined PDE system for $z(u,v)$ yields \eqref{Liouville-gen-soln}, where
\[ \phi(u) = \int e^{\rho(u)}\, du, \qquad \psi(v) = \int \tfrac{1}{2} e^{-\sigma(v)}\, dv. \]
\end{remark}

By making the change of variables $\tilde{u} = \phi(u)$, $\tilde{v} = \psi(v)$, we can write \eqref{apply-Liouville-gen-soln} as
\[ f(u,v) + g(u,v) + 2h(u) = \ln \left( \frac{2}{(u+v)^2} \right). \]
Thus we have
\begin{align}
f + g & = \ln \left( \frac{2}{(u+v)^2} \right) - 2 h(u) \notag \\
& = \ln \left( \frac{2}{(H(u))^2(u+v)^2} \right). \label{fg-expr}
\end{align}
Now we can write equation \eqref{PDE2} as
\begin{align*}
H'(u) & = 2 f_v e^{2f} e^{-(f+g)} \\
& = \frac{(e^{2f})_v (H(u))^2}{(f+g)_{uv}} \text{\ (from equation \eqref{PDE3})};
\end{align*}
therefore,
\[ \frac{H'(u)}{(H(u))^2} (f+g)_{uv} = (e^{2f})_v. \]
Integrating with respect to $v$ yields
\[ \frac{H'(u)}{(H(u))^2} (f+g)_{u} = (e^{2f}) + k(u) \]
for some function $k(u)$.  But from equation \eqref{fg-expr}, we have
\[ (f+g)_u = -2\left( \frac{H'(u)}{H(u)} + \frac{1}{u+v} \right). \]
Therefore, 
\begin{align*}
e^{2f} & = \frac{H'(u)}{(H(u))^2} (f+g)_{u} - k(u) \\
& =  \frac{-2H'(u)}{(H(u))^2} \left( \frac{H'(u)}{H(u)} + \frac{1}{u+v} \right) - k(u) \\
& = \frac{-2H'(u)}{(H(u))^2 (u+v)} - \tilde{k}(u), 
\end{align*}
where
\[ \tilde{k}(u) = k(u) + 2\frac{(H'(u))^2}{(H(u))^3}. \]
Thus the Maurer-Cartan forms are given by \eqref{MC-forms-1}, \eqref{MC-forms-2}, \eqref{MC-forms-3}, where
\begin{gather}
H = H(u) \neq 0, \notag \\
f = \tfrac{1}{2} \ln \left(\frac{-2H'(u)}{(H(u))^2 (u+v)} - k(u) \right), \label{fgH-case2} \\
g = \ln \left( \frac{2}{(H(u))^2(u+v)^2} \right) - \tfrac{1}{2} \ln \left(\frac{-2H'(u)}{(H(u))^2 (u+v)} - k(u) \right), \notag
\end{gather}
and $k(u)$ is an arbitrary function of $u$.  Thus solutions of the PDE system \eqref{PDE1}-\eqref{PDE3} with $H(u) \neq 0$ depend on two arbitrary functions of one variable---i.e., the functions $h(u), k(u)$.

\section{Explicit parametrizations}\label{main-result-sec}

In this section we use the Maurer-Cartan forms obtained in \S \ref{fundamental-forms-sec} to derive explicit local parametrizations for all totally quasi-umbilic timelike surfaces.  The Maurer-Cartan equations \eqref{MC-forms-def} can be written in matrix form as
\begin{equation}
 \begin{bmatrix} d\bx & d\f_1 & d\f_2 & d\f_3 \end{bmatrix} = \begin{bmatrix} \bx & \f_1 & \f_2 & \f_3 \end{bmatrix} \begin{bmatrix} 0 & 0 & 0 & 0 \\[0.1in] \eta^1 & \eta^1_1 & \eta^1_2 & \eta^1_3 \\[0.1in] \eta^2 & \eta^2_1 & \eta^2_2 & \eta^2_3 \\[0.1in] \eta^3 & \eta^3_1 & \eta^3_2 & \eta^3_3 \end{bmatrix}. \label{MC-matrix-eqn}
\end{equation}

\subsection{Case 1: $H(u) \equiv 0$}  In this case, equation \eqref{MC-matrix-eqn} takes the form
\begin{equation*}
 \begin{bmatrix} d\bx & d\f_1 & d\f_2 & d\f_3 \end{bmatrix} = \begin{bmatrix} \bx & \f_1 & \f_2 & \f_3 \end{bmatrix} \begin{bmatrix} 0 & 0 & 0 & 0 \\[0.1in] e^{f(u)}\, du & -f'(u)\, du & 0 & 0 \\[0.1in] e^{-f(u)}\, dv & 0 & f'(u)\, du & -e^{f(u)}\, du \\[0.1in] 0 & -e^{f(u)}\, du & 0 & 0 \end{bmatrix}. 
\end{equation*}
This is equivalent to the compatible, overdetermined PDE system
\begin{gather}
 \begin{bmatrix} \bx & \f_1 & \f_2 & \f_3 \end{bmatrix}_u  = \begin{bmatrix} \bx & \f_1 & \f_2 & \f_3 \end{bmatrix} \begin{bmatrix} 0 & 0 & 0 & 0 \\[0.1in] e^{f(u)} & -f'(u) & 0 & 0 \\[0.1in] 0 & 0 & f'(u) & -e^{f(u)} \\[0.1in] 0 & -e^{f(u)} & 0 & 0 \end{bmatrix}, \label{case1-u-pde} \\[0.2in]
 \begin{bmatrix} \bx & \f_1 & \f_2 & \f_3 \end{bmatrix}_v  = \begin{bmatrix} \bx & \f_1 & \f_2 & \f_3 \end{bmatrix} \begin{bmatrix} 0 & 0 & 0 & 0 \\[0.1in] 0 & 0 & 0 & 0 \\[0.1in] e^{-f(u)} & 0 & 0 & 0 \\[0.1in] 0 & 0 & 0 & 0 \end{bmatrix}  \label{case1-v-pde}
\end{gather}
for the vector-valued functions $\bx, \f_1, \f_2, \f_3$ on $\Sigma$.
These equations can be integrated explicitly, and the general solution for $\bx(u,v)$ is
\[
\bx(u,v) = \bx^0 +  u e^{f(0)} \f^0_1 + \left[ v + \tfrac{1}{2} \int_0^u \left( \int_0^{\sigma} e^{2f(\tau)} d\tau \right)^2 d\sigma \right] e^{-f(0)}\f^0_2 
- \left[ \int_0^u \left( \int_0^{\sigma} e^{2f(\tau)} d\tau \right) d\sigma \right] \f^0_3, 
\]
where $\bx^0, \f_1^0, \f_2^0, \f_3^0$ represent the initial conditions
\[ \begin{bmatrix} \bx(0,0) & \f_1(0,0) & \f_2(0,0) & \f_3(0,0) \end{bmatrix} = \begin{bmatrix} \bx^0 & \f^0_1 & \f^0_2 & \f^0_3 \end{bmatrix}. \]
Thus $\Sigma$ is a cylindrical surface of the form
\begin{equation}
 \bx(u,v) = \alpha(u) + v \f^0_2, \label{cylinder-surface}
\end{equation}
where the $v$-coordinate curves are null lines parallel to the null vector $\f^0_2$, and the generating curve $\alpha(u)$ is the null curve
\[
 \alpha(u) = \bx^0 +  u e^{f(0)} \f^0_1 + \left[ \tfrac{1}{2} \int_0^u \left( \int_0^{\sigma} e^{2f(\tau)} d\tau \right)^2 d\sigma \right]  e^{-f(0)}\f^0_2 - \left[ \int_0^u \left( \int_0^{\sigma} e^{2f(\tau)} d\tau \right) d\sigma \right] \f^0_3 .
\]

We can give a slightly simpler parametrization for $\Sigma$ (not by null coordinates) as \eqref{cylinder-surface} with
\[ \alpha(u) = \bx^0 + u e^{f(0)} \f^0_1 - \left[ \int_0^u \left( \int_0^{\sigma} e^{2f(\tau)} d\tau \right) d\sigma \right] \f^0_3. \]
Then $\alpha$ may be viewed as the graph in the $\f_1^0 \f_3^0$-plane of an arbitrary function $F(u)$ with $F''(u) < 0$.  We summarize this result, noting that
\begin{itemize}
\item the $\f_1^0 \f_3^0$-plane is a lightlike plane, 
\item $\f_3^0$ spans the unique direction in this plane which is orthogonal to $\f_2^0$, and 
\item $\alpha'(u)$ is always linearly independent from $\f_3^0$,
\end{itemize}
 as:

\begin{proposition}\label{H=0-prop}
Let $\Sigma$ be a totally quasi-umbilic, regular timelike surface in $\R^{1,2}$ with mean curvature $H \equiv 0$.  Then $\Sigma$ is a cylinder over a convex curve $\alpha$ which is contained in a lightlike plane, with null rulings transverse to the plane of $\alpha$; moreover, the tangent vector to $\alpha$ is nowhere orthogonal to the rulings.  Conversely, any such cylinder is totally quasi-umbilic with mean curvature $H \equiv 0$.
\end{proposition}

\begin{remark}
Since the surfaces of Proposition \ref{H=0-prop} have $K \equiv 0$ as well, this proposition shows that there exists a family of non-planar surfaces in $\R^{1,2}$, parametrized by one arbitrary function of one variable, all of which have Gauss and mean curvature identically equal to zero!
\end{remark}

\subsection{Case 2: $H(u) \neq 0$} In this case, the PDE system represented by equation \eqref{MC-matrix-eqn} is considerably more complicated and generally cannot be completely solved analytically, but there is still enough information in \eqref{MC-matrix-eqn} to give a geometric description of the solution surfaces.  The $v$-component of the corresponding PDE system is
\begin{equation}
\begin{bmatrix} \bx & \f_1 & \f_2 & \f_3 \end{bmatrix}_v  = \begin{bmatrix} \bx & \f_1 & \f_2 & \f_3 \end{bmatrix} \begin{bmatrix} 0 & 0 & 0 & 0 \\[0.1in] 0 & -f_v & 0 & 0 \\[0.1in] e^g & 0 & f_v & -He^g \\[0.1in] 0 & -He^g & 0 & 0 \end{bmatrix},  \label{case2-v-pde}
\end{equation}
where $f, g, H$ are given by \eqref{fgH-case2}.  In particular, we have
\[ (\f_2)_v = f_v \f_2, \]
and so 
\[ \f_2(u,v) = e^{f} \bar{\f}_2(u) \]
for some null vector field $\bar{\f}_2(u)$ along $\Sigma$ which is a function of $u$ alone.  This means that along any $v$-coordinate curve $\bx(u_0, v)$ in $\Sigma$, the vectors $\f_2(u_0, v)$ are all parallel.  Moreover, we have
\[ \bx_v(u,v) = e^g \f_2(u,v) = e^{f+g} \bar{\f}_2(u); \]
it follows that the $v$-coordinate curve $\bx(u_0, v)$ is a null line parallel to $\bar{\f}_2(u_0)$.  So if we let $\alpha$ be the null curve 
\[ \alpha(u) = \bx(u, 0) \]
and define $\bar{\f}_2(u)$ to be the null vector field 
\[ \bar{\f}_2(u) = \f_2(u, 0) \]
along $\alpha$, then we can reparametrize $\Sigma$ (not necessarily by null coordinates) as
\[ \bx(u,v) = \alpha(u) + v\bar{\f}_2(u). \]
Therefore, any such surface is ruled, and the rulings are all null lines.  Moreover, the $u$-component of the PDE system corresponding to \eqref{MC-matrix-eqn} 
implies that:
\begin{itemize}
\item The vectors $\alpha'(u)$ and $\alpha''(u)$ are linearly independent for all $u$; we call such a null curve {\em nondegenerate}. 
\item The vectors $\bar{\f}_2(u)$ and $\bar{\f}'_2(u)$ are linearly independent if and only if $H(u) \neq 0$.
\end{itemize}

Conversely, suppose that $\Sigma$ is a timelike ruled surface whose rulings are all null lines.  Let $\alpha$ be a null curve contained in $\Sigma$ which is transverse to the rulings, and suppose further that $\alpha$ is nondegenerate.  Define vector fields $\bar{\f}_1(u), \bar{\f}_2(u), \bar{\f}_3(u)$ along $\alpha$ by the conditions that:
\begin{itemize}
\item $\bar{\f}_1(u) = \alpha'(u)$;
\item $\bar{\f}_2(u)$ spans the ruling of $\Sigma$ passing through $\alpha(u)$, and $\langle \bar{\f}_1(u), \bar{\f}_2(u) \rangle = 1$;
\item $\bar{\f}_3(u)$ is a unit normal vector to $\Sigma$ at $\alpha(u)$.
\end{itemize}
Then we can parametrize $\Sigma$ as
\[ \bx(u,v) = \alpha(u) + v \bar{\f}_2(u), \]
and there exists an adapted null frame field $(\f_1, \f_2, \f_3)$ along $\Sigma$ with the property that
\[ \f_1(u,0) = \bar{\f}_1(u), \qquad \f_2(u,0) = \bar{\f}_2(u), \qquad \f_3(u,0) = \bar{\f}_3(u). \]

The $u$-component of the Maurer-Cartan equation \eqref{MC-matrix-eqn} for the adapted null frame field $(\f_1, \f_2, \f_3)$ evaluated along $\alpha$ implies that there exist functions $h_{ij}(u)$ such that
\[ \begin{bmatrix} \bar{\f}'_1(u) & \bar{\f}'_2(u) & \bar{\f}'_3(u) \end{bmatrix} = \begin{bmatrix} \bar{\f}_1(u) & \bar{\f}_2(u) & \bar{\f}_3(u) \end{bmatrix} \begin{bmatrix} h_{11}(u) & 0 & h_{32}(u) \\[0.1in] 0 & -h_{11}(u) & h_{31}(u) \\[0.1in] h_{31}(u) & h_{32}(u) & 0 \end{bmatrix}, \]
and the condition that $\alpha$ is nondegenerate implies that $h_{31}(u) \neq 0$ for all $u$.

In order to show that $\Sigma$ is totally quasi-umbilic, we will compute the shape operator $S_{\bx}$ at any point $\bx \in \Sigma$ with respect to the basis $(\bx_u, \bx_v)$ for $T_{\bx}\Sigma$.  First we compute:
\begin{align*}
\bx_u & = \bar{\f}_1(u) - v h_{11}(u) \bar{\f}_2(u) + v h_{32}(u) \bar{\f}_3(u), \\
\bx_v & = \bar{\f}_2(u).
\end{align*}
It is straightforward to check that the unit normal vector field $\f_3$ to $\Sigma$ is given by
\[ \f_3(u,v) = \bar{\f}_3(u) + v h_{32}(u) \bar{\f}_2(u). \]
So the action of the shape operator on the basis $(\bx_u, \bx_v)$ is given by:
\begin{align*}
S_{\bx}(\bx_u) & = -(\f_3)_u \\
   & = -\bar{\f}_3'(u) - v(h'_{32}(u) \bar{\f}_2(u) + h_{32}(u) \bar{\f}'_2(u)) \\
   & = -h_{32}(u) \bar{\f}_1(u) + (-h_{31}(u) - v h'_{32}(u) + v h_{32}(u) h_{11}(u)) \bar{\f}_2(u) - v (h_{32}(u))^2 \bar{\f}_3(u)       \\
   & =  -h_{32}(u) ( \bar{\f}_1(u) - v h_{11}(u) \bar{\f}_2(u) + v h_{32}(u) \bar{\f}_3(u)) - (h_{31}(u) + v h'_{32}(u)) \bar{\f}_2(u)     \\
   & = -h_{32}(u) \bx_u - (h_{31}(u) + v h'_{32}(u)) \bx_v, \\
S_{\bx}(\bx_v) & = -(\f_3)_v \\
& = -h_{32}(u) \bar{\f}_2(u) \\
& = -h_{32}(u) \bx_v.
\end{align*}
Therefore, the matrix of $S_{\bx}$ with respect to the basis $(\bx_u, \bx_v)$ is
\[ \begin{bmatrix} -h_{32}(u) & 0 \\[0.1in] -h_{31}(u) - v h'_{32}(u) & -h_{32}(u) \end{bmatrix} . \]
Since $h_{31}(u) \neq 0$, it follows that $\Sigma$ is quasi-umbilic almost everywhere, and if $h'_{32}(u) \neq 0$, then $\Sigma$ contains a curve of umbilic points defined by the equation $v = -\frac{h_{31}(u)}{h'_{32}(u)}$.

Together with Proposition \ref{H=0-prop}, this proves the following classification theorem:

\begin{theorem}\label{main-thm}
Let $\Sigma$ be a totally quasi-umbilic, regular timelike surface in $\R^{1,2}$. Then $\Sigma$ is a ruled surface whose rulings are all null lines, with the additional property that any null curve $\alpha$ in $\Sigma$ which is transverse to the rulings is nondegenerate.  Conversely, given any nondegenerate null curve $\alpha(u)$ in $\R^{1,2}$ and any null vector field $\bar{\f}_2(u)$ along $\alpha$ which is linearly independent from $\alpha'(u)$ for all $u$, the immersed ruled surface
\[ \bx(u,v) = \alpha(u) + v \bar{\f}_2(u) \]
is totally quasi-umbilic, possibly containing a curve of umbilic points.  
\end{theorem}

\section{examples}\label{examples-sec}

In light of Theorem \ref{main-thm}, any totally quasi-umbilic, regular timelike surface $\Sigma$ has a parametrization of the form
\[ \bx(u,v) = \alpha(u) + v \bar{\f}_2(u), \]
where $\alpha$ is a nondegenerate null curve whose tangent vector $\bar{\f}_1(u) = \alpha'(u)$ is a null vector which is linearly independent from the null vector $\bar{\f}_2(u)$.  The nondegeneracy of $\alpha$ implies that $\bar{\f}_1'(u)$ must be linearly independent from $\bar{\f}_1(u)$.

For convenience, we can reparametrize the curve $\alpha(u)$ and scale the vector field $\bar{\f}_2(u)$ so as to arrange that both $\bar{\f}_1(u)$ and $\bar{\f}_2(u)$  have {\em Euclidean} length 2; i.e., we can assume that
\[ \bar{\f}_1(u) = (1, \cos \theta_1(u), \sin \theta_1(u)), \qquad 
\bar{\f}_2(u) = (1, \cos \theta_2(u), \sin \theta_2(u)) \]
for some functions $\theta_1(u), \theta_2(u)$.  (These vectors will no longer necessarily satisfy the condition that $\langle \bar{\f}_1(u), \bar{\f}_2(u) \rangle = 1$, but this condition is not important at this point.)  Then the linear independence conditions of the previous paragraph are simply that $\theta_1'(u) \neq 0$ and $\theta_1(u) \neq \theta_2(u)$ for all $u$.  Moreover, $\Sigma$ has Gauss and mean curvature $H \equiv K \equiv 0$ if and only if $\theta_2'(u) \equiv 0$.

\subsection{Case 1: $H(u) \equiv 0$}  In this case $\theta_2(u)$ must be constant, while $\theta_1(u)$ must be either strictly increasing or strictly decreasing and never equal to $\theta_2(u)$.

\begin{example}\label{example-1}
Let 
\[   \theta_1(u) = 2\tan^{-1} u , \qquad  \theta_2(u) = \pi . \] 
Then we have $\theta_1'(u) > 0$ for all $u$, and since $-\pi < \theta_1(u) < \pi$, the linear independence requirements are satisfied.  The corresponding null vector fields are given by
\begin{align*}
\bar{\f}_1(u) & = \left(1, \cos \left(2 \tan^{-1} u \right), \sin \left(2 \tan^{-1} u \right) \right) \\
& = \left(1, \frac{1 - u^2}{1 + u^2}, \frac{2u}{1 + u^2} \right)
\end{align*}
and
$\bar{\f}_2(u) = (1, -1, 0). $
Thus we have 
\[ \alpha(u) = \int \bar{\f}_1(u)\, du = \left(u,\, 2 \tan^{-1} u - u,\, \ln (1 + u^2) \right), \]
and the corresponding surface $\Sigma$ is given by
\begin{align*} 
\bx(u,v) & = \alpha(u) + v \bar{\f}_2(u) \\
& = \left(u + v, \, 2 \tan^{-1} u - u - v, \, \ln (1 + u^2) \right).
\end{align*}
Note that we could rewrite $\bx(u,v)$ as
\begin{align*}
\bx(u,v) & = \left( \tan^{-1} u, \, \tan^{-1} u, \ln(1 + u^2) \right) +  (u - \tan^{-1} u + v, -(u - \tan^{-1} u + v), 0) \\
& = \left( \tan^{-1} u, \, \tan^{-1} u, \ln(1 + u^2) \right) + \tilde{v} (1, -1, 0),
\end{align*}
where $\tilde{v} = u - \tan^{-1} u + v.$  This realizes $\Sigma$ as a cylinder over the convex curve
\[ \alpha(u) = (\tan^{-1} u) (1, 1, 0) + \ln(1 + u^2) (0,0,1) \]
in the lightlike plane spanned by the vectors $(1,1,0)$ and $(0,0,1)$.  Making the substitution $\sigma = \tan^{-1} u$ shows that this curve is the graph of the function $F(\sigma) = \ln(\sec^2 \sigma)$.  (See Figure \ref{example-1-fig}; note that the $x^0$-axis is drawn as the vertical axis.)
\begin{figure}[h]
\includegraphics[width=3in]{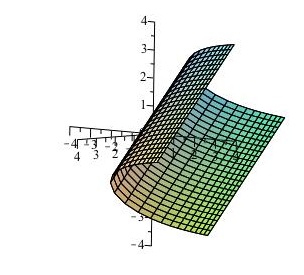}
\caption{The surface of Example \ref{example-1}}
\label{example-1-fig}
\end{figure}

\end{example}

\subsection{Case 2: $H(u) \neq 0$} In this case, $\theta_1(u)$ and $\theta_2(u)$ must both be either strictly increasing or strictly decreasing, and we must have $\theta_1(u) \neq \theta_2(u)$ for all $u$.

\begin{example}\label{example-2}
Let
\[ \theta_1(u) = u, \qquad \theta_2(u) = u + \pi. \]
The corresponding null vector fields are given by
\[ \bar{\f}_1(u) = (1,\, \cos u,\, \sin u), \qquad \bar{\f}_2(u) = (1,\, -\cos u,\, -\sin u). \]
Thus we have
\[ \alpha(u) = \int \bar{\f}_1(u)\, du = \left(u,\, \sin u,\, -\cos u \right), \]
and the corresponding surface $\Sigma$ is given by
\begin{align*} 
\bx(u,v) & = \alpha(u) + v \bar{\f}_2(u) \\
& = \left(u + v, \, \sin u - v \cos u, \, -\cos u - v \sin u \right).
\end{align*}
(See Figure \ref{example-2-fig} for embedded and immersed portions of this surface; note that the $x^0$-axis is drawn as the vertical axis.)
\begin{figure}[h]
\includegraphics[width=2.5in]{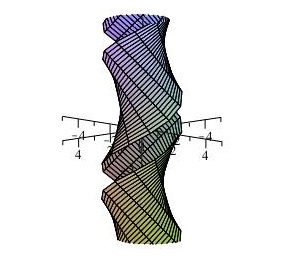}
\includegraphics[width=2.5in]{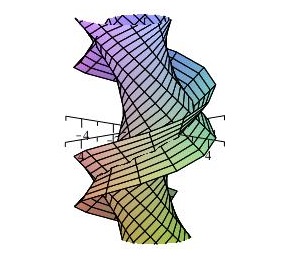}
\caption{Embedded and immersed portions of the surface of Example \ref{example-2}}
\label{example-2-fig}
\end{figure}

In order to compute shape operator $S_{\bx}$, first we compute
\begin{align*}
\bx_u & = (1, \, \cos u + v \sin u, \, \sin u - v \cos u) \\
\bx_v & = (1, \, -\cos u,\, -\sin u).
\end{align*}
The unit normal vector to $\Sigma$ is given by
\begin{align*}
 \f_3(u,v) & = \frac{\bx_u \times \bx_v}{\sqrt{-\langle \bx_u \times \bx_v, \bx_u \times \bx_v\rangle}} \\
 &  =  \left(-\frac{v}{2}, \, -\sin u + \frac{v}{2} \cos u, \, \cos u + \frac{v}{2} \sin u \right).
\end{align*}
(Note that the cross product here denotes the Minkowski cross product \eqref{Minkowski-cross-product}.)  Now the action of the shape operator on the basis $(\bx_u, \bx_v)$ for $T_{\bx}\Sigma$ is given by:
\begin{align*}
S_{\bx}(\bx_u) & = -(\f_3)_u \\
   & = \left( 0, \, \cos u + \tfrac{1}{2} \sin u, \, \sin u - \tfrac{1}{2} \cos u \right) \\
   & = \tfrac{1}{2} \bx_u - \tfrac{1}{2} \bx_v, \\
S_{\bx}(\bx_v) & = -(\f_3)_v \\
& = \left( \tfrac{1}{2},\,  -\tfrac{1}{2} \cos u, \, \tfrac{1}{2} \sin u \right) \\
& = \tfrac{1}{2} \bx_v.
\end{align*}
Therefore the matrix of $S_{\bx}$ with respect to the basis $(\bx_u, \bx_v)$ for $T_{\bx}\Sigma$ is
\[ \begin{bmatrix} \tfrac{1}{2} & 0 \\[0.1in] -\tfrac{1}{2} & \tfrac{1}{2} \end{bmatrix} . \]
Thus we see that $\Sigma$ has constant Gauss and mean curvatures $K \equiv \tfrac{1}{4}$, $H \equiv \tfrac{1}{2}$.  Moreover, $\Sigma$ contains no umbilic points.
\end{example}

\begin{example}\label{example-3}
Let
\[ \theta_1(u) = u^3 + u, \qquad \theta_2(u) = u^3 + u + \pi. \]
The corresponding null vector fields are given by
\[ \bar{\f}_1(u) = (1,\, \cos (u^3 + u),\, \sin (u^3 + u)), \qquad \bar{\f}_2(u) = (1,\, -\cos (u^3 + u),\, -\sin (u^3 + u)). \]
Define
\[ \alpha(u) = \int_0^u \bar{\f}_1(s)\, ds = \left(u,\, \int_0^u \cos (s^3 + s)\, ds,\, \int_0^u \sin (s^3 + s)\, ds \right), \]
and the corresponding surface $\Sigma$ is given by
\begin{align*} 
\bx(u,v) & = \alpha(u) + v \bar{\f}_2(u) \\
& =  \left(u + v,\, \int_0^u \cos (s^3 + s)\, ds - v\cos(u^3 + u),\, \int_0^u \sin (s^3 + s)\, ds - v\sin(u^3 + u) \right).
\end{align*}

In order to compute shape operator $S_{\bx}$, first we compute
\begin{align*}
\bx_u & = (1, \, \cos (u^3 + u) + v (3 u^2 + 1) \sin (u^3 + u), \, \sin (u^3 + u) - v (3 u^2 + 1) \cos (u^3 + u)) \\
\bx_v & = (1, \, -\cos (u^3 + u),\, -\sin (u^3 + u)).
\end{align*}
The unit normal vector to $\Sigma$ is given by
\begin{align*}
 \f_3(u,v) & = \frac{\bx_u \times \bx_v}{\sqrt{-\langle \bx_u \times \bx_v, \bx_u \times \bx_v\rangle}} \\
& =  \left( 0, \, -\sin (u^3 + u), \, \cos (u^3 + u) \right) + \frac{v}{2} (3u^2 + 1) \left( -1, \, \cos (u^3 + u), \, \sin (u^3 + u) \right).
\end{align*}
A computation shows that the action of the shape operator on the basis $(\bx_u, \bx_v)$ for $T_{\bx}\Sigma$ is given by:
\begin{align*}
S_{\bx}(\bx_u) 
   & = \tfrac{1}{2} (3u^2 + 1) \bx_u  + \left( 3uv - \tfrac{1}{2} (3u^2 + 1) \right)  \bx_v, \\[0.1in]
S_{\bx}(\bx_v) 
& = \tfrac{1}{2} (3u^2 + 1) \bx_v.
\end{align*}
Therefore the matrix of $S_{\bx}$ with respect to the basis $(\bx_u, \bx_v)$ for $T_{\bx}\Sigma$ is
\[ \begin{bmatrix} \tfrac{1}{2} (3u^2 + 1) & 0 \\[0.1in] 3uv - \tfrac{1}{2} (3u^2 + 1) & \tfrac{1}{2} (3u^2 + 1) \end{bmatrix} . \]
Thus we see that $\Sigma$ has Gauss and mean curvatures 
\[ K = \tfrac{1}{4} (3u^2 + 1)^2, \qquad H = \tfrac{1}{2} (3u^2 + 1). \]
Moreover, $\Sigma$ contains a curve of umbilic points defined by the equation
\[ 3uv - \tfrac{1}{2} (3u^2 + 1) = 0. \]
See Figure \ref{example-3-fig} for various views of embedded and immersed portions of this surface, including the curve of umbilic points; note that the $x^0$-axis is drawn as the vertical axis.
\begin{figure}[h]
\includegraphics[width=2.5in]{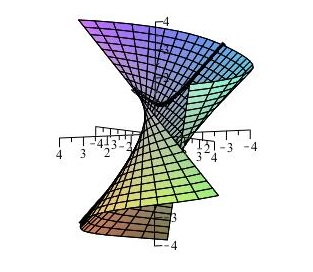}
\includegraphics[width=2.5in]{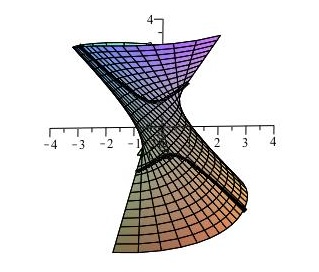}
\includegraphics[width=2.5in]{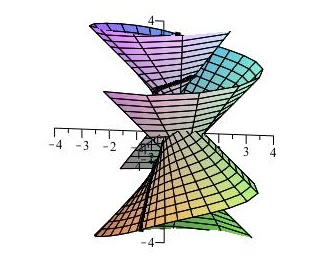}
\includegraphics[width=2.5in]{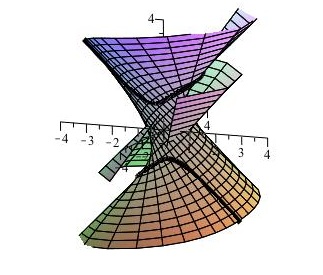}
\caption{Embedded and immersed portions of the surface of Example \ref{example-3}}
\label{example-3-fig}
\end{figure}

\end{example}

\begin{example}\label{example-4}
The previous two examples were immersed surfaces; we can obtain an embedded example by perturbing the totally umbilic hyperboloid $\langle \bv, \bv \rangle = 1$.  We begin with the standard parametrization
\[ \bx(u,v) = (0,\, \cos u ,\, \sin u) + v (1,\, -\sin u , \, \cos u) \]
for the hyperboloid.  (Unlike our previous examples, the generating curve $\alpha(u) = (0,\, \cos u ,\, \sin u)$ is not a null curve, but this is not crucial.)  We will modify this example by replacing the parameter $u$ in the ruling vector field by the function $u + \tfrac{1}{2} \sin u$.

\begin{remark}
Any small perturbation of the parameter $u$ would suffice, as long as the perturbation is $2\pi$-periodic and has a derivative of magnitude less than 1.  Alternatively, we could perturb the parameter $u$ in the generating curve, or perform some combination of the two perturbations.
\end{remark}

Thus we will consider the surface $\Sigma$ with parametrization
\[ \bx(u,v) = (0,\, \cos u ,\, \sin u) + v \left( 1,\, -\sin (u + \tfrac{1}{2} \sin u) , \, \cos (u + \tfrac{1}{2} \sin u) \right). \]
Computations similar to those of the previous example show that $\Sigma$ is totally quasi-umbilic, with a curve of umbilic points defined by the equation
\begin{multline*}
 [8 \sin(u - \sin u) + 8 \sin u - \sin(2u + \sin u) + \sin (2u - \sin u) - 2\sin(\sin u)] v \\
   +  [ 8\cos(u + \tfrac{1}{2} \sin u) + 8 \cos(u - \tfrac{1}{2} \sin u)] = 0. \end{multline*}
See Figure \ref{example-4-fig} for two views of this surface, including the curve of umbilic points (which has two components); note that the $x^0$-axis is drawn as the vertical axis.
\begin{figure}[h]
\includegraphics[width=2.5in]{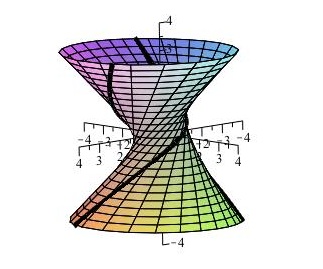}
\includegraphics[width=2.5in]{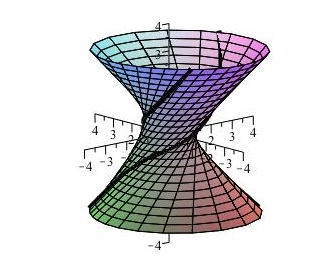}
\caption{The surface of Example \ref{example-4}}
\label{example-4-fig}
\end{figure}

\end{example}

\bibliographystyle{amsplain}
\bibliography{quasi-umbilic-bib}

\providecommand{\bysame}{\leavevmode\hbox to3em{\hrulefill}\thinspace}
\providecommand{\MR}{\relax\ifhmode\unskip\space\fi MR }
\providecommand{\MRhref}[2]{%
  \href{http://www.ams.org/mathscinet-getitem?mr=#1}{#2}
}
\providecommand{\href}[2]{#2}
\begin{thebibliography}{1}

\bibitem{CI09}
Jeanne~N. Clelland and Thomas~A. Ivey, \emph{B\"acklund transformations and
  {D}arboux integrability for nonlinear wave equations}, Asian J. Math.
  \textbf{13} (2009), no.~1, 15--64.

\bibitem{doCarmo76}
Manfredo~P. do~Carmo, \emph{Differential geometry of curves and surfaces},
  Prentice-Hall Inc., Englewood Cliffs, N.J., 1976, Translated from the
  Portuguese.

\bibitem{IL03}
Thomas~A. Ivey and J.~M. Landsberg, \emph{Cartan for beginners: differential
  geometry via moving frames and exterior differential systems}, Graduate
  Studies in Mathematics, vol.~61, American Mathematical Society, Providence,
  RI, 2003.

\bibitem{Magid05}
Martin~A. Magid, \emph{Lorentzian isothermic surfaces in {${\bf R}^n_j$}},
  Rocky Mountain J. Math. \textbf{35} (2005), no.~2, 627--640.

\bibitem{Magid09}
\bysame, \emph{Blaschke's problem for timelike surfaces}, Results Math.
  \textbf{56} (2009), no.~1-4, 335--355.

\bibitem{Milnor83}
Tilla~Klotz Milnor, \emph{Harmonic maps and classical surface theory in
  {M}inkowski {$3$}-space}, Trans. Amer. Math. Soc. \textbf{280} (1983), no.~1,
  161--185.

\end{thebibliography}

\end{document}